\documentclass[11pt]{amsart}
\usepackage{amsmath,amssymb}
\newtheorem{theorem}{Theorem}[section]
\newtheorem{proposition}[theorem]{Proposition}
\newtheorem{corollary}[theorem]{Corollary}
\newtheorem{lemma}[theorem]{Lemma}

\begin{document}

\title{A boundary value problem for minimal Lagrangian graphs}
\author{Simon Brendle and Micah Warren}
%\begin{abstract}
%Let $\Omega$ and $\tilde{\Omega}$ be uniformly convex domains in $\mathbb{R}^n$ with smooth boundary. We show that there exists a diffeomorphism $f: \Omega \to \tilde{\Omega}$ such that the graph $\Sigma = \{(x,f(x)): x \in \Omega\}$ is a minimal Lagrangian submanifold of $\mathbb{R}^n \times \mathbb{R}^n$.
%\end{abstract}

\maketitle 

\section{Introduction}

Consider the product $\mathbb{R}^n \times \mathbb{R}^n$ equipped with the Euclidean metric. The product $\mathbb{R}^n \times \mathbb{R}^n$ has a natural complex structure, which is given by 
\[J \frac{\partial}{\partial x_k} = \frac{\partial}{\partial y_k}, \quad J \frac{\partial}{\partial y_k} = -\frac{\partial}{\partial x_k}.\] 
The associated symplectic structure is given by 
\[\omega = \sum_{k=1}^n dx_k \wedge dy_k.\] 
A submanifold $\Sigma \subset \mathbb{R}^n \times \mathbb{R}^n$ is called Lagrangian if $\omega|_\Sigma = 0$.

In this paper, we study a boundary value problem for minimal Lagrangian graphs in $\mathbb{R}^n \times \mathbb{R}^n$. To that end, we fix two domains $\Omega,\tilde{\Omega} \subset \mathbb{R}^n$ with smooth boundary. Given a diffeomorphism $f: \Omega \to \tilde{\Omega}$, we consider its graph $\Sigma = \{(x,f(x)): x \in \Omega\} \subset \mathbb{R}^n \times \mathbb{R}^n$. We consider the problem of finding a diffeomorphism $f: \Omega \to \tilde{\Omega}$ such that $\Sigma$ is Lagrangian and has zero mean curvature. Our main result asserts that such a map exists if $\Omega$ and $\tilde{\Omega}$ are uniformly convex:

\begin{theorem} 
\label{existence}
Let $\Omega$ and $\tilde{\Omega}$ be uniformly convex domains in $\mathbb{R}^n$ with smooth boundary. Then there exists a diffeomorphism $f: \Omega \to \tilde{\Omega}$ such that the graph 
\[\Sigma = \{(x,f(x)): x \in \Omega\}\] 
is a minimal Lagrangian submanifold of $\mathbb{R}^n \times \mathbb{R}^n$.
\end{theorem}

Minimal Lagrangian submanifolds were first studied by Harvey and Lawson \cite{Harvey-Lawson}, and have attracted considerable interest in recent years. Yuan \cite{Yuan} has proved a Bernstein-type theorem for minimal Lagrangian graphs over $\mathbb{R}^n$. A similar result was established by Tsui and Wang \cite{Tsui-Wang}. Smoczyk and Wang have used the mean curvature flow to deform certain Lagrangian submanifolds to minimal Lagrangian submanifolds (see \cite{Smoczyk}, \cite{Smoczyk-Wang}, \cite{Wang}). In \cite{Brendle}, the first author studied a boundary value problem for minimal Lagrangian graphs in $\mathbb{H}^2 \times \mathbb{H}^2$, where $\mathbb{H}^2$ denotes the hyperbolic plane.

In order to prove Theorem \ref{existence}, we reduce the problem to the solvability of a fully nonlinear PDE. As above, we assume that $\Omega$ and $\tilde{\Omega}$ are uniformly convex domains in $\mathbb{R}^n$ with smooth boundary. Moreover, suppose that $f$ is a diffeomorphism from $\Omega$ to $\tilde{\Omega}$. The graph $\Sigma = \{(x,f(x)): x \in \Omega\}$ is Lagrangian if and only if there exists a function $u: \Omega \to \mathbb{R}$ such that $f(x) = \nabla u(x)$. In that case, the Lagrangian angle of $\Sigma$ is given by $F(D^2 u(x))$. Here, $F$ is a real-valued function on the space of symmetric $n \times n$ matrices which is defined as follows: if $M$ is a symmetric $n \times n$ matrix, then $F(M)$ is defined by 
\[F(M) = \sum_{k=1}^n \arctan(\lambda_k),\] 
where $\lambda_1, \hdots, \lambda_n$ denote the eigenvalues of $M$. 

By a result of Harvey and Lawson (see \cite{Harvey-Lawson}, Proposition 2.17), $\Sigma$ has zero mean curvature if and only if the Lagrangian angle is constant; that is, 
\begin{equation} 
\label{pde}
F(D^2 u(x)) = c 
\end{equation}
for all $x \in \Omega$. Hence, we are led to the following problem: \\

$(\star)$ \textit{Find a convex function $u: \Omega \to \mathbb{R}$ and a constant $c \in (0,\frac{n\pi}{2})$ such that $\nabla u$ is a diffeomorphism from $\Omega$ to $\tilde{\Omega}$ and $F(D^2 u(x)) = c$ for all $x \in \Omega$.} \\

Caffarelli, Nirenberg, and Spruck \cite{Caffarelli-Nirenberg-Spruck} have obtained an existence result for solutions of (\ref{pde}) under Dirichlet boundary conditions. In this paper, we study a different boundary condition, which is analogous to the second boundary value problem for the Monge-Amp\`ere equation. 

In dimension $2$, P.~Delano\"e \cite{Delanoe} proved that the second boundary value problem for the Monge-Amp\`ere equation has a unique smooth solution, provided that both domains are uniformly convex. This result was generalized to higher dimensions by L.~Caffarelli \cite{Caffarelli} and J.~Urbas \cite{Urbas1}. In 2001, J.~Urbas \cite{Urbas2} described a general class of Hessian equations for which the second boundary value problem admits a unique smooth solution.

In Section 2, we establish a-priori estimates for solutions of $(\star)$. In Section 3, we prove that all solutions of $(\star)$ are non-degenerate (that is, the linearized operator is invertible). In Section 4, we use the continuity method to show that $(\star)$ has at least one solution. From this, Theorem \ref{existence} follows. Finally, in Section 5, we prove a uniqueness result for $(\star)$.

The first author is grateful to Professors Philippe Delano\"e and John Urbas for discussions. The first author was partially supported by the National Science Foundation under grants DMS-0605223 and DMS-0905628. The second author was partially supported by a Liftoff Fellowship from the Clay Mathematics Institute.

\section{A priori estimates for solutions of $(\star)$}

In this section, we prove a-priori estimates for solutions of $(\star)$. 

Let $\Omega$ and $\tilde{\Omega}$ be uniformly convex domains in $\mathbb{R}^n$ with smooth boundary. Moreover, suppose that $u$ is a convex function such that $\nabla u$ is a diffeomorphism from $\Omega$ to $\tilde{\Omega}$ and $F(D^2 u(x))$ is constant. For each point $x \in \Omega$, we define a symmetric $n \times n$-matrix $A(x) = \{a_{ij}(x): 1 \leq i,j \leq n\}$ by 
\[A(x) = \big [ I + (D^2 u(x))^2 \big ]^{-1}.\] 
Clearly, $A(x)$ is positive definite for all $x \in \Omega$. 

\begin{lemma}
\label{upper.bound.for.Lagrangian.angle}
We have 
\[\frac{n\pi}{2} - F(D^2 u(x)) \geq \arctan \bigg ( \frac{\text{\rm vol}(\Omega)^{1/n}}{\text{\rm vol}(\tilde{\Omega})^{1/n}} \bigg )\] 
for all points $x \in \Omega$.
\end{lemma}

\textbf{Proof.} 
Since $\nabla u$ is a diffeomorphism from $\Omega$ to $\tilde{\Omega}$, we have 
\[\int_\Omega \det D^2 u(x) \, dx = \text{\rm vol}(\tilde{\Omega}).\] 
Therefore, we can find a point $x_0 \in \Omega$ such that 
\[\det D^2 u(x_0) \leq \frac{\text{\rm vol}(\tilde{\Omega})}{\text{\rm vol}(\Omega)}.\] 
Hence, if we denote by $\lambda_1 \leq \lambda_2 \leq \hdots \leq \lambda_n$ the eigenvalues of $D^2 u(x_0)$, then we have 
\[\lambda_1 \leq \frac{\text{\rm vol}(\tilde{\Omega})^{1/n}}{\text{\rm vol}(\Omega)^{1/n}}.\] 
This implies
\begin{align*} 
\frac{n\pi}{2} - F(D^2 u(x_0)) &= \sum_{k=1}^n \arctan \Big ( \frac{1}{\lambda_k} \Big ) \\ 
&\geq \arctan \Big ( \frac{1}{\lambda_1} \Big ) \\ 
&\geq \arctan \bigg ( \frac{\text{\rm vol}(\Omega)^{1/n}}{\text{\rm vol}(\tilde{\Omega})^{1/n}} \bigg ). 
\end{align*}
Since $F(D^2 u(x))$ is constant, the assertion follows. \\

\begin{lemma} 
\label{bound.for.smallest.eigenvalue}
Let $x$ be an arbitrary point in $\Omega$, and let $\lambda_1 \leq \lambda_2 \leq \hdots \leq \lambda_n$ be the eigenvalues of $D^2 u(x)$. Then 
\[\frac{1}{\lambda_1} \geq \tan \bigg [ \frac{1}{n} \, \arctan \bigg ( \frac{\text{\rm vol}(\Omega)^{1/n}}{\text{\rm vol}(\tilde{\Omega})^{1/n}} \bigg ) \bigg ].\] 
\end{lemma}

\textbf{Proof.} 
Using Lemma \ref{upper.bound.for.Lagrangian.angle}, we obtain 
\begin{align*} 
n \, \arctan \Big ( \frac{1}{\lambda_1} \Big ) 
&\geq \sum_{k=1}^n \arctan \Big ( \frac{1}{\lambda_k} \Big ) \\ 
&= \frac{n\pi}{2} - F(D^2 u(x)) \\ 
&\geq \arctan \bigg ( \frac{\text{\rm vol}(\Omega)^{1/n}}{\text{\rm vol}(\tilde{\Omega})^{1/n}} \bigg ). 
\end{align*}
From this, the assertion follows easily. \\

By Proposition \ref{boundary.defining.function}, we can find a smooth function $h: \Omega \to \mathbb{R}$ such that $h(x) = 0$ for all $x \in \partial \Omega$ and 
\begin{equation} 
\label{convexity.of.h}
\sum_{i,j=1}^n \partial_i \partial_j h(x) \, w_i \, w_j \geq \theta \, |w|^2 
\end{equation}
for all $x \in \Omega$ and all $w \in \mathbb{R}^n$. Similarly, there exists a smooth function $\tilde{h}: \tilde{\Omega} \to \mathbb{R}$ such that $\tilde{h}(y) = 0$ for all $y \in \partial \tilde{\Omega}$ and 
\begin{equation} 
\label{convexity.of.h.tilde}
\sum_{i,j=1}^n \partial_i \partial_j \tilde{h}(y) \, w_i \, w_j \geq \theta \, |w|^2 
\end{equation}
for all $y \in \tilde{\Omega}$ and all $w \in \mathbb{R}^n$. For abbreviation, we choose a positive constant $C_1$ such that 
\[C_1 \, \theta \, \sin^2 \bigg [ \frac{1}{n} \, \arctan \bigg ( \frac{\text{\rm vol}(\Omega)^{1/n}}{\text{\rm vol}(\tilde{\Omega})^{1/n}} \bigg ) \bigg ] = 1.\] 
We then have the following estimate:

\begin{lemma} 
\label{pde.for.h}
We have 
\[\sum_{i,j=1}^n a_{ij}(x) \, \partial_i \partial_j h(x) \geq \frac{1}{C_1}\] 
for all $x \in \Omega$.
\end{lemma}

\textbf{Proof.} 
Fix a point $x_0 \in \Omega$, and let $\lambda_1 \leq \lambda_2 \leq \hdots \leq \lambda_n$ be the eigenvalues of $D^2 u(x_0)$. It follows from (\ref{convexity.of.h}) that 
\[\sum_{i,j=1}^n a_{ij}(x_0) \, \partial_i \partial_j h(x_0) \geq \theta \, \sum_{k=1}^n \frac{1}{1 + \lambda_k^2} \geq \theta \, \frac{1}{1 + \lambda_1^2}.\] 
Using Lemma \ref{bound.for.smallest.eigenvalue}, we obtain 
\[\frac{1}{1 + \lambda_1^2} \geq \sin^2 \bigg [ \frac{1}{n} \, \arctan \bigg ( \frac{\text{\rm vol}(\Omega)^{1/n}}{\text{\rm vol}(\tilde{\Omega})^{1/n}} \bigg ) \bigg ] = \frac{1}{C_1 \, \theta}.\] 
Putting these facts together, the assertion follows. \\

In the next step, we differentiate the identity $F(D^2 u(x)) = \text{\rm constant}$ with respect to $x$. To that end, we need the following well-known fact: 

\begin{lemma} 
\label{variation.of.F}
Let $M(t)$ be a smooth one-parameter family of symmetric $n \times n$ matrices. Then 
\[\frac{d}{dt} F(M(t)) \Big |_{t=0} = \text{\rm tr} \big [ (I + M(0)^2)^{-1} \, M'(0) \big ].\] 
Moreover, if $M(0)$ is positive definite, then we have 
\[\frac{d^2}{dt^2} F(M(t)) \Big |_{t=0} \leq \text{\rm tr} \big [ (I + M(0)^2)^{-1} \, M''(0) \big ].\] 
\end{lemma}

\begin{proposition} 
We have 
\begin{equation} 
\label{first.derivative.of.F}
\sum_{i,j=1}^n a_{ij}(x) \, \partial_i \partial_j \partial_k u(x) = 0 
\end{equation}
for all $x \in \Omega$. Moreover, we have 
\begin{equation} 
\label{second.derivative.of.F}
\sum_{i,j,k,l=1}^n a_{ij}(x) \, \partial_i \partial_j \partial_k \partial_l u(x) \, w_k \, w_l \geq 0 
\end{equation}
for all $x \in \Omega$ and all $w \in \mathbb{R}^n$.
\end{proposition}

\textbf{Proof.} 
Fix a point $x_0 \in \Omega$ and a vector $w \in \mathbb{R}^n$. It follows from Lemma \ref{variation.of.F} that 
\[0 = \frac{d}{dt} F \big ( D^2 u(x_0 + tw) \big ) \Big |_{t=0} = \sum_{i,j,k=1}^n a_{ij}(x) \, \partial_i \partial_j \partial_k u(x_0) \, w_k.\] 
Moreover, since the matrix $D^2 u(x_0)$ is positive definite, we have 
\[0 = \frac{d^2}{dt^2} F \big ( D^2 u(x_0 + tw) \big ) \Big |_{t=0} \leq \sum_{i,j,k,l=1}^n a_{ij}(x) \, \partial_i \partial_j \partial_k \partial_l u(x_0) \, w_k \, w_l.\] From this, the assertion follows. \\

\begin{proposition} 
\label{general.estimate}
Fix a smooth function $\Phi: \Omega \times \tilde{\Omega} \to \mathbb{R}$, and define $\varphi(x) = \Phi(x,\nabla u(x))$. Then 
\[\bigg | \sum_{i,j=1}^n a_{ij}(x) \, \partial_i \partial_j \varphi(x) \bigg | \leq C\] 
for all $x \in \Omega$. Here, $C$ is a positive constant that depends only on the second order partial derivatives of $\Phi$.
\end{proposition}

\textbf{Proof.} 
The partial derivatives of the function $\varphi(x)$ are given by 
\[\partial_i \varphi(x) = \sum_{k=1}^n \Big ( \frac{\partial}{\partial y_k} \Phi \Big )(x,\nabla u(x)) \, \partial_i \partial_k u(x) + \Big ( \frac{\partial}{\partial x_i} \Phi \Big ) (x,\nabla u(x)).\] 
This implies 
\begin{align*} 
\partial_i \partial_j \varphi(x) 
&= \sum_{k=1}^n \Big ( \frac{\partial}{\partial y_k} \Phi \Big ) (x,\nabla u(x)) \, \partial_i \partial_j \partial_k u(x) \\ 
&+ \sum_{k,l=1}^n \Big ( \frac{\partial^2}{\partial y_k \partial y_l} \Phi \Big ) (x,\nabla u(x)) \, \partial_i \partial_k u(x) \, \partial_j \partial_l u(x) \\ 
&+ \sum_{k=1}^n \Big ( \frac{\partial^2}{\partial x_j \partial y_k} \Phi \Big ) (x,\nabla u(x)) \, \partial_i \partial_k u(x) \\ 
&+ \sum_{l=1}^n \Big ( \frac{\partial^2}{\partial x_i \partial y_l} \Phi \Big ) (x,\nabla u(x)) \, \partial_j \partial_l u(x) \\ 
&+ \Big ( \frac{\partial^2}{\partial x_i \partial x_j} \Phi \Big ) (x,\nabla u(x)). 
\end{align*} 
Using (\ref{first.derivative.of.F}), we obtain 
\begin{align*} 
&\sum_{i,j=1}^n a_{ij}(x) \, \partial_i \partial_j \varphi(x) \\ 
&= \sum_{i,j,k,l=1}^n a_{ij}(x) \, \Big ( \frac{\partial^2}{\partial y_k \partial y_l} \Phi \Big ) (x,\nabla u(x)) \, \partial_i \partial_k u(x) \, \partial_j \partial_l u(x) \\ 
&+ 2 \sum_{i,j,k=1}^n a_{ij}(x) \, \Big ( \frac{\partial^2}{\partial x_j \partial y_k} \Phi \Big ) (x,\nabla u(x)) \, \partial_i \partial_k u(x) \\ 
&+ \sum_{i,j=1}^n a_{ij}(x) \, \Big ( \frac{\partial^2}{\partial x_i \partial x_j} \Phi \Big ) (x,\nabla u(x)). 
\end{align*} 
We now fix a point $x_0 \in \Omega$. Without loss of generality, we may assume that $D^2 u(x_0)$ is a diagonal matrix. This implies 
\begin{align*} 
\sum_{i,j=1}^n a_{ij}(x_0) \, \partial_i \partial_j \varphi(x_0) 
&= \sum_{k=1}^n \frac{\lambda_k^2}{1+\lambda_k^2} \, \Big ( \frac{\partial^2}{\partial y_k^2} \Phi \Big ) (x_0,\nabla u(x_0)) \\ 
&+ 2 \sum_{k=1}^n \frac{\lambda_k}{1+\lambda_k^2} \, \Big ( \frac{\partial^2}{\partial x_k \partial y_k} \Phi \Big ) (x_0,\nabla u(x_0)) \\ 
&+ \sum_{k=1}^n \frac{1}{1+\lambda_k^2} \, \Big ( \frac{\partial^2}{\partial x_k^2} \Phi \Big ) (x_0,\nabla u(x_0)), 
\end{align*} 
where $\lambda_k = \partial_k \partial_k u(x_0)$. Thus, we conclude that 
\[\bigg | \sum_{i,j=1}^n a_{ij}(x_0) \, \partial_i \partial_j \varphi(x_0) \bigg | \leq C,\] 
as claimed. \\

We next consider the function $H(x) = \tilde{h}(\nabla u(x))$. The following estimate is an immediate consequence of Proposition \ref{general.estimate}:

\begin{corollary} 
\label{pde.for.tilde.h}
There exists a positive constant $C_2$ such that 
\[\bigg | \sum_{i,j=1}^n a_{ij}(x) \, \partial_i \partial_j H(x) \bigg | \leq C_2\] 
for all $x \in \Omega$.
\end{corollary}

\begin{proposition}
\label{bound.for.H}
We have $H(x) \geq C_1C_2 \, h(x)$ for all $x \in \Omega$.
\end{proposition}

\textbf{Proof.} 
Using Lemma \ref{pde.for.h} and Corollary \ref{pde.for.tilde.h}, we obtain 
\[\sum_{i,j=1}^n a_{ij}(x) \, \partial_i \partial_j(H(x) - C_1C_2 \, h(x)) \leq 0\] 
for all $x \in \Omega$. Hence, the function $H(x) - C_1C_2 \, h(x)$ attains its minimum on $\partial \Omega$. Thus, we conclude that $H(x) - C_1C_2 \, h(x) \geq 0$ for all $x \in \Omega$. \\

\begin{corollary}
\label{bound.for.gradient.of.H}
We have 
\[\langle \nabla h(x),\nabla H(x) \rangle \leq C_1C_2 \, |\nabla h(x)|^2\] 
for all $x \in \partial \Omega$.
\end{corollary}

\begin{proposition}
\label{normal.derivative}
Fix a smooth function $\Phi: \Omega \times \tilde{\Omega} \to \mathbb{R}$, and define $\varphi(x) = \Phi(x,\nabla u(x))$. Then 
\[|\langle \nabla \varphi(x),\nabla \tilde{h}(\nabla u(x)) \rangle| \leq C\] 
for all $x \in \partial \Omega$. Here, $C$ is a positive constant that depends only on $C_1,C_2$, and the first order partial derivatives of $\Phi$.
\end{proposition}

\textbf{Proof.} 
A straightforward calculation yields 
\begin{align*} 
\langle \nabla \varphi(x),\nabla \tilde{h}(\nabla u(x)) \rangle 
&= \sum_{k=1}^n \Big ( \frac{\partial}{\partial x_k} \Phi \Big ) (x,\nabla u(x)) \, (\partial_k \tilde{h})(\nabla u(x)) \\
&+ \sum_{k=1}^n \Big ( \frac{\partial}{\partial y_k} \Phi \Big ) (x,\nabla u(x)) \, \partial_k H(x) 
\end{align*}
for all $x \in \Omega$. By Corollary \ref{bound.for.gradient.of.H}, we have $|\nabla H(x)| \leq C_1 C_2 \, |\nabla h(x)|$ for all points $x \in \partial \Omega$. Putting these facts together, the assertion follows. \\

\begin{proposition}
\label{obliqueness}
We have 
\begin{align*} 
0 
&< \sum_{k,l=1}^n \partial_k \partial_l u(x) \, (\partial_k \tilde{h})(\nabla u(x)) \, (\partial_l \tilde{h})(\nabla u(x)) \\ 
&\leq C_1C_2 \, \langle \nabla h(x),\nabla \tilde{h}(\nabla u(x)) \rangle 
\end{align*} 
for all $x \in \partial \Omega$.
\end{proposition}

\textbf{Proof.} 
Note that the function $H$ vanishes along $\partial \Omega$ and is negative in the interior of $\Omega$. Hence, for each point $x \in \partial \Omega$, the vector $\nabla H(x)$ is a positive multiple of $\nabla h(x)$. Since $u$ is convex, we obtain 
\begin{align*} 
0 
&< \sum_{k,l=1}^n \partial_k \partial_l u(x) \, (\partial_k \tilde{h})(\nabla u(x)) \, (\partial_l \tilde{h})(\nabla u(x)) \\ 
&= \langle \nabla H(x),\nabla \tilde{h}(\nabla u(x)) \rangle \\ 
&= \frac{\langle \nabla h(x),\nabla H(x) \rangle}{|\nabla h(x)|^2} \, \langle \nabla h(x),\nabla \tilde{h}(\nabla u(x)) \rangle 
\end{align*} 
for all $x \in \partial \Omega$. In particular, we have $\langle h(x),\nabla \tilde{h}(\nabla u(x)) \rangle > 0$ for all points $x \in \partial \Omega$. The assertion follows now from Corollary \ref{bound.for.gradient.of.H}. \\

\begin{proposition} 
\label{uniform.obliqueness}
There exists a positive constant $C_4$ such that 
\[\langle \nabla h(x),\nabla \tilde{h}(\nabla u(x)) \rangle 
\geq \frac{1}{C_4}\] 
for all $x \in \partial \Omega$.
\end{proposition}

\textbf{Proof.} 
We define a function $\chi(x)$ by 
\[\chi(x) = \langle \nabla h(x),\nabla \tilde{h}(\nabla u(x)) \rangle.\] 
By Proposition \ref{general.estimate}, we can find a positive constant $C_3$ such that 
\[\bigg | \sum_{i,j=1}^n a_{ij}(x) \, \partial_i \partial_j \chi(x) \bigg | \leq C_3\] 
for all $x \in \Omega$. Using Lemma \ref{pde.for.h}, we obtain  
\[\sum_{i,j=1}^n a_{ij}(x) \, \partial_i \partial_j (\chi(x) - C_1 C_3 \, h(x)) \leq 0\] 
for all $x \in \Omega$. Hence, there exists a point $x_0 \in \partial \Omega$ such that 
\[\inf_{x \in \Omega} (\chi(x) - C_1 C_3 \, h(x)) = \inf_{x \in \partial \Omega} \chi(x) = \chi(x_0).\] 
It follows from Proposition \ref{obliqueness} that $\chi(x_0) > 0$. Moreover, we can find a nonnegative real number $\mu$ such that 
\[\nabla \chi(x_0) = (C_1 C_3 - \mu) \, \nabla h(x_0).\] 
A straightforward calculation yields 
\begin{align} 
\label{gradient.chi}
\langle \nabla \chi(x),\nabla \tilde{h}(\nabla u(x)) \rangle 
&= \sum_{i,j=1}^n \partial_i \partial_j h(x) \, (\partial_i \tilde{h})(\nabla u(x)) \, (\partial_j \tilde{h})(\nabla u(x)) \notag \\ 
&+ \sum_{i,j=1}^n (\partial_i \partial_j \tilde{h})(\nabla u(x)) \, \partial_i h(x) \, \partial_j H(x) 
\end{align} 
for all $x \in \partial \Omega$. Using (\ref{convexity.of.h}), we obtain 
\[\sum_{i,j=1}^n \partial_i \partial_j h(x) \, (\partial_i \tilde{h})(\nabla u(x)) \, (\partial_j \tilde{h})(\nabla u(x)) \geq \theta \, |\nabla \tilde{h}(\nabla u(x))|^2\] 
for all $x \in \partial \Omega$. Since $\nabla H(x)$ is a positive multiple of $\nabla h(x)$, we have 
\[\sum_{i,j=1}^n (\partial_i \partial_j \tilde{h})(\nabla u(x)) \, \partial_i h(x) \, \partial_j H(x) \geq 0\] 
for all $x \in \partial \Omega$. Substituting these inequalities into (\ref{gradient.chi}) gives 
\[\langle \nabla \chi(x),\nabla \tilde{h}(\nabla u(x)) \rangle \geq \theta \, |\nabla \tilde{h}(\nabla u(x))|^2\] 
for all $x \in \partial \Omega$. From this, we deduce that 
\begin{align*} 
(C_1C_3 - \mu) \, \chi(x_0) 
&= (C_1C_3 - \mu) \, \langle \nabla h(x_0),\nabla \tilde{h}(\nabla u(x_0)) \rangle \\ 
&= \langle \nabla \chi(x_0),\nabla \tilde{h}(\nabla u(x_0)) \rangle \\ 
&\geq \theta \, |\nabla \tilde{h}(\nabla u(x_0))|^2. 
\end{align*} 
Since $\mu \geq 0$ and $\chi(x_0) > 0$, we conclude that 
\[\chi(x_0) \geq \frac{\theta}{C_1C_3} \, |\nabla \tilde{h}(\nabla u(x_0))|^2 \geq \frac{1}{C_4}\] 
for some positive constant $C_4$. This completes the proof of Proposition \ref{uniform.obliqueness}. \\

\begin{lemma}
\label{M}
Suppose that 
\[\sum_{k,l=1}^n \partial_k \partial_l u(x) \, w_k \, w_l \leq M \, |w|^2\] 
for all $x \in \partial \Omega$ and all $w \in T_x(\partial \Omega)$. Then 
\begin{align*} 
\sum_{k,l=1}^n \partial_k \partial_l u(x) \, w_k \, w_l 
&\leq M \, \bigg | w - \frac{\langle \nabla h(x),w \rangle}{\langle \nabla h(x),\nabla \tilde{h}(\nabla u(x)) \rangle} \, \nabla \tilde{h}(\nabla u(x)) \bigg |^2 \\ 
&+ C_1 C_2 C_4 \, \langle \nabla h(x),w \rangle^2 
\end{align*} 
for all $x \in \partial \Omega$ and all $w \in \mathbb{R}^n$.
\end{lemma}

\textbf{Proof.} 
Fix a point $x \in \partial \Omega$ and a vector $w \in \mathbb{R}^n$. Morever, let 
\[z = w - \frac{\langle \nabla h(x),w \rangle}{\langle \nabla h(x),\nabla \tilde{h}(\nabla u(x)) \rangle} \, \nabla \tilde{h}(\nabla u(x)).\] 
Clearly, $\langle \nabla h(x),z \rangle = 0$; hence $z \in T_x(\partial \Omega)$. This implies 
\[\sum_{k,l=1}^n \partial_k \partial_l u(x) \, (\partial_k \tilde{h})(\nabla u(x)) \, z_l = \langle \nabla H(x),z \rangle = 0.\] From this we deduce that 
\begin{align*} 
&\sum_{k,l=1}^n \partial_k \partial_l u(x) \, w_k \, w_l - \sum_{k,l=1}^n \partial_k \partial_l u(x) \, z_k \, z_l \\ 
&= \frac{\langle \nabla h(x),w \rangle^2}{\langle \nabla h(x),\nabla \tilde{h}(\nabla u(x)) \rangle^2} \, \sum_{k,l=1}^n \partial_k \partial_l u(x) \, (\partial_k \tilde{h})(\nabla u(x)) \, (\partial_l \tilde{h})(\nabla u(x)). 
\end{align*} 
It follows from Proposition \ref{obliqueness} and Proposition \ref{uniform.obliqueness} that 
\begin{align*} 
&\frac{\langle \nabla h(x),w \rangle^2}{\langle \nabla h(x),\nabla \tilde{h}(\nabla u(x)) \rangle^2} \, \sum_{k,l=1}^n \partial_k \partial_l u(x) \, (\partial_k \tilde{h})(\nabla u(x)) \, (\partial_l \tilde{h})(\nabla u(x)) \\ 
&\leq C_1 C_2 \, \frac{\langle \nabla h(x),w \rangle^2}{\langle \nabla h(x),\nabla \tilde{h}(\nabla u(x)) \rangle} 
\leq C_1 C_2 C_4 \, \langle \nabla h(x),w \rangle^2. 
\end{align*} 
Moreover, we have 
\[\sum_{k,l=1}^n \partial_k \partial_l u(x) \, z_k \, z_l \leq M \, |z|^2\] 
by definition of $M$. Putting these facts together, the assertion follows. \\

\begin{proposition} 
\label{estimate.for.tangential.derivatives}
There exists a positive constant $C_9$ such that 
\[\sum_{k,l=1}^n \partial_k \partial_l u(x) \, w_k \, w_l \leq C_9 \, |w|^2\] 
for all $x \in \partial \Omega$ and all $w \in T_x(\partial \Omega)$.
\end{proposition} 

\textbf{Proof.} 
Let 
\[M = \sup \bigg \{ \sum_{k,l=1}^n \partial_k \partial_l u(x) \, z_k \, z_l: x \in \partial \Omega, \, z \in T_x (\partial \Omega), \, |z| = 1 \bigg \}.\] 
By compactness, we can find a point $x_0 \in \partial \Omega$ and a unit vector $w \in T_{x_0}(\partial \Omega)$ such that 
\[\sum_{k,l=1}^n \partial_k \partial_l u(x_0) \, w_k \, w_l = M.\] 
We define a function $\psi: \Omega \to \mathbb{R}$ by 
\[\psi(x) = \sum_{k,l=1}^n \partial_k \partial_l u(x) \, w_k \, w_l\] 
for all $x \in \Omega$. Moreover, we define functions $\varphi_1: \Omega \to \mathbb{R}$ and $\varphi_2: \Omega \to \mathbb{R}$ by 
\[\varphi_1(x) = \bigg | w - \frac{\langle \nabla h(x),w \rangle}{\eta(\langle \nabla h(x),\nabla \tilde{h}(\nabla u(x)) \rangle)} \, \nabla \tilde{h}(\nabla u(x)) \bigg |^2\] 
and 
\[\varphi_2(x) = \langle \nabla h(x),w \rangle^2\] 
for all $x \in \Omega$. Here, $\eta: \mathbb{R} \to \mathbb{R}$ is a smooth cutoff function satisfying $\eta(s) = s$ for $s \geq \frac{1}{C_4}$ and $\eta(s) \geq \frac{1}{2C_4}$ for all $s \in \mathbb{R}$. 

The inequality (\ref{second.derivative.of.F}) implies that 
\[\sum_{i,j=1}^n a_{ij}(x) \, \partial_i \partial_j \psi(x) \geq 0\] 
for all $x \in \Omega$. Moreover, by Proposition \ref{general.estimate}, there exists a positive constant $C_5$ such that 
\[\bigg | \sum_{i,j=1}^n a_{ij}(x) \, \partial_i \partial_j \varphi_1(x) \bigg | \leq C_5\] 
and 
\[\bigg | \sum_{i,j=1}^n a_{ij}(x) \, \partial_i \partial_j \varphi_2(x) \bigg | \leq C_5\] 
for all $x \in \Omega$. Hence, the function 
\begin{align*} 
g(x) &= \psi(x) - M \, \varphi_1(x) - C_1 C_2 C_4 \, \varphi_2(x) \\ 
&+ C_1C_5 \, (M + C_1 C_2 C_4) \, h(x) 
\end{align*} 
satisfies
\begin{equation} 
\label{pde.for.g}
\sum_{i,j=1}^n a_{ij}(x) \, \partial_i \partial_j g(x) \geq 0 
\end{equation}
for all $x \in \Omega$.

It follows from Proposition \ref{uniform.obliqueness} that 
\[\varphi_1(x) = \bigg | w - \frac{\langle \nabla h(x),w \rangle}{\langle \nabla h(x),\nabla \tilde{h}(\nabla u(x)) \rangle} \, \nabla \tilde{h}(\nabla u(x)) \bigg |^2\] 
for all $x \in \partial \Omega$. Using Lemma \ref{M}, we obtain 
\[\psi(x) \leq M \, \varphi_1(x) + C_1 C_2 C_4 \, \varphi_2(x)\] 
for all $x \in \partial \Omega$. Therefore, we have $g(x) \leq 0$ for all $x \in \partial \Omega$. Using the inequality (\ref{pde.for.g}) and the maximum principle, we conclude that $g(x) \leq 0$ for all $x \in \Omega$. 

On the other hand, we have $\varphi_1(x_0) = 1$, $\varphi_2(x_0) = 0$, and $\psi(x_0) = M$. From this, we deduce that $g(x_0) = 0$. Therefore, the function $g$ attains its global maximum at the point $x_0$. This implies $\nabla g(x_0) = \mu \, \nabla h(x_0)$ for some nonnegative real number $\mu$. From this, we deduce that 
\begin{equation} 
\label{estimate.1}
\langle \nabla g(x_0),\nabla \tilde{h}(\nabla u(x_0)) \rangle = \mu \, \langle \nabla h(x_0),\nabla \tilde{h}(\nabla u(x_0)) \rangle \geq 0. 
\end{equation} 
By Proposition \ref{normal.derivative}, we can find a positive constant $C_6$ such that 
\[|\langle \nabla \varphi_1(x),\nabla \tilde{h}(\nabla u(x)) \rangle| \leq C_6\] 
for all $x \in \partial \Omega$. Hence, we can find positive constants $C_7$ and $C_8$ such that
\begin{align} 
\label{estimate.2}
\langle \nabla g(x),\nabla \tilde{h}(\nabla u(x)) \rangle 
&= \langle \nabla \psi(x),\nabla \tilde{h}(\nabla u(x)) \rangle \notag \\ 
&- M \, \langle \nabla \varphi_1(x),\nabla \tilde{h}(\nabla u(x)) \rangle \notag \\ 
&- C_1 C_2 C_4 \, \langle \nabla \varphi_2(x),\nabla \tilde{h}(\nabla u(x)) \rangle \\ 
&+ C_1C_5 \, (M + C_1 C_2 C_4) \, \langle \nabla h(x),\nabla \tilde{h}(\nabla u(x)) \rangle \notag \\ 
&\leq \langle \nabla \psi(x),\nabla \tilde{h}(\nabla u(x)) \rangle + C_7 \, M + C_8 \notag 
\end{align} 
for all $x \in \partial \Omega$. Combining (\ref{estimate.1}) and (\ref{estimate.2}), we conclude that 
\begin{equation} 
\label{estimate.5}
\langle \nabla \psi(x_0),\nabla \tilde{h}(\nabla u(x_0)) \rangle + C_7 \, M + C_8 \geq 0. 
\end{equation}

A straightforward calculation shows that 
\begin{align} 
\label{estimate.6}
&\sum_{k,l=1}^n \partial_k \partial_l H(x_0) \, w_k \, w_l \notag \\ 
&= \sum_{i,k,l=1}^n (\partial_i \tilde{h})(\nabla u(x_0)) \, \partial_i \partial_k \partial_l u(x_0) \, w_k \, w_l \\ 
&+ \sum_{i,j,k,l=1}^n (\partial_i \partial_j \tilde{h})(\nabla u(x_0)) \, \partial_i \partial_k u(x_0) \, \partial_j \partial_l u(x_0) \, w_k \, w_l. \notag
\end{align} 
Since $H$ vanishes along $\partial \Omega$, we have 
\[\sum_{k,l=1}^n \partial_k \partial_l H(x_0) \, w_k \, w_l = -\langle \nabla H(x_0),I\!I(w,w) \rangle,\] 
where $I\!I(\cdot,\cdot)$ denotes the second fundamental form of $\partial \Omega$ at $x_0$. 
Using the estimate $|\nabla H(x_0)| \leq C_1C_2 \, |\nabla h(x_0)|$, we obtain  
\[\sum_{k,l=1}^n \partial_k \partial_l H(x_0) \, w_k \, w_l \leq C_1 C_2 \, |\nabla h(x_0)| \, |I\!I(w,w)|.\] 
Moreover, we have 
\[\sum_{i,k,l=1}^n (\partial_i \tilde{h})(\nabla u(x_0)) \, 
\partial_i \partial_k \partial_l u(x_0) \, w_k \, w_l = \langle \nabla \psi(x_0),\nabla \tilde{h}(\nabla u(x_0)) \rangle.\] 
Finally, it follows from (\ref{convexity.of.h.tilde}) that 
\begin{align*} 
&\sum_{i,j,k,l=1}^n (\partial_i \partial_j \tilde{h})(\nabla u(x_0)) \, \partial_i \partial_k u(x_0) \, \partial_j \partial_l u(x_0) \, w_k \, w_l \\ 
&\geq \theta \, \sum_{i,j,k,l=1}^n \partial_i \partial_k u(x_0) \, \partial_j \partial_l u(x_0) \, w_i \, w_j \, w_k \, w_l = \theta \, M^2. 
\end{align*} 
Substituting these inequalities into (\ref{estimate.6}), we obtain 
\begin{align*} 
C_1C_2 \, |\nabla h(x_0)| \, |I\!I(w,w)| 
&\geq \sum_{k,l=1}^n \partial_k \partial_l H(x_0) \, w_k \, w_l \\ 
&\geq \langle \nabla \psi(x_0),\nabla \tilde{h}(\nabla u(x_0)) \rangle + \theta \, M^2 \\ 
&\geq \theta \, M^2 - C_7 \, M - C_8. 
\end{align*} 
Therefore, we have $M \leq C_9$ for some positive constant $C_9$. This completes the proof of Proposition \ref{estimate.for.tangential.derivatives}. \\

\begin{corollary} 
\label{boundary.c.2.estimate}
There exists a positive constant $C_{10}$ such that 
\[\sum_{k,l=1}^n \partial_k \partial_l u(x) \, w_k \, w_l \leq C_{10} \, |w|^2\] 
for all $x \in \partial \Omega$ and all $w \in \mathbb{R}^n$.
\end{corollary}

\textbf{Proof.} 
It follows from Lemma \ref{M} that 
\begin{align*} 
\sum_{k,l=1}^n \partial_k \partial_l u(x) \, w_k \, w_l 
&\leq C_9 \, \bigg | w - \frac{\langle \nabla h(x),w \rangle}{\langle \nabla h(x),\nabla \tilde{h}(\nabla u(x)) \rangle} \, \nabla \tilde{h}(\nabla u(x)) \bigg |^2 \\ 
&+ C_1 C_2 C_4 \, \langle \nabla h(x),w \rangle^2
\end{align*} 
for all $x \in \partial \Omega$ and all $w \in \mathbb{R}^n$. Hence, the assertion follows from Proposition \ref{uniform.obliqueness}. \\

The interior $C^2$ estimate follows from Corollary \ref{boundary.c.2.estimate} and (\ref{second.derivative.of.F}):

\begin{proposition} 
\label{interior.c.2.estimate}
We have 
\[\sum_{k,l=1}^n \partial_k \partial_l u(x) \, w_k \, w_l \leq C_{10} \, |w|^2\] 
for all $x \in \Omega$ and all $w \in \mathbb{R}^n$.
\end{proposition}

\textbf{Proof.} 
Fix a unit vector $w \in \mathbb{R}^n$, and define 
\[\psi(x) = \sum_{k,l=1}^n \partial_k \partial_l u(x) \, w_k \, w_l.\] 
The inequality (\ref{second.derivative.of.F}) implies that 
\[\sum_{i,j=1}^n a_{ij}(x) \, \partial_i \partial_j \psi(x) \geq 0\] 
for all $x \in \Omega$. Using the maximum principle, we obtain 
\[\sup_{x \in \Omega} \psi(x) = \sup_{x \in \partial \Omega} \psi(x) \leq C_{10}.\] 
This completes the proof. \\

Once we have a uniform $C^2$ bound, we can show that $u$ is uniformly convex:

\begin{corollary}
\label{uniform.convexity}
There exists a positive constant $C_{11}$ such that 
\[\sum_{k,l=1}^n \partial_k \partial_l u(x) \, w_k \, w_l \geq \frac{1}{C_{11}} \, |w|^2\] 
for all $x \in \Omega$ and all $w \in \mathbb{R}^n$.
\end{corollary}

\textbf{Proof.} 
By assumption, the map $f(x) = \nabla u(x)$ is a diffeomorphism from $\Omega$ to $\tilde{\Omega}$. Let $g: \tilde{\Omega} \to \Omega$ denote the inverse of $f$. Then $Dg(y) = \big [ Df(x) \big ]^{-1}$, where $x = g(y)$. Since the matrix $Df(x) = D^2 u(x)$ is positive definite for all $x \in \Omega$, we conclude that the matrix $Dg(y)$ is positive definite for all $y \in \tilde{\Omega}$. Hence, there exists a convex function $v: \tilde{\Omega} \to \mathbb{R}$ such that $g(y) = \nabla v(y)$. The function $v$ satisfies $F(D^2 v(y)) = \frac{n\pi}{2} - F(D^2 u(x))$, where $x = g(y)$. Since $F(D^2 u(x))$ is constant, it follows that $F(D^2 v(y))$ is constant. By Proposition \ref{interior.c.2.estimate}, the eigenvalues of $D^2 v(y)$ are uniformly bounded from above. From this, the assertion follows. \\

In the next step, we show that the second derivatives of $u$ are uniformly bounded in $C^\gamma(\overline{\Omega})$. To that end, we use results of G.~Lieberman and N.~Trudinger \cite{Lieberman-Trudinger}. In the remainder of this section, we describe how the problem $(\star)$ can be rewritten so as to fit into the framework of Lieberman and Trudinger.

We begin by choosing a smooth cutoff function $\eta: \mathbb{R} \to [0,1]$ such that 
\[\begin{cases} \eta(s) = 0 & \text{\rm for $s \leq 0$} \\ 
\eta(s) = 1 & \text{\rm for $\frac{1}{C_{11}} \leq s \leq C_{10}$} \\ 
\eta(s) = 0 & \text{\rm for $s \geq 2C_{10}$}. \end{cases}\] 
There exists a unique function $\psi: \mathbb{R} \to \mathbb{R}$ satisfying $\psi(1) = \frac{\pi}{4}$, $\psi'(1) = \frac{1}{2}$, and $\psi''(s) = -\frac{2s}{(1+s^2)^2} \, \eta(s) \leq 0$ for all $s \in \mathbb{R}$. Clearly, $\psi(s) = \arctan(s)$ for $\frac{1}{C_{11}} \leq s \leq C_{10}$. Moreover, it is easy to see that $\frac{1}{1 + 4C_{10}^2} \leq \psi'(s) \leq 1$ for all $s \in \mathbb{R}$. If $M$ is a symmetric $n \times n$ matrix, we define 
\[\Psi(M) = \sum_{k=1}^n \psi(\lambda_k),\] 
where $\lambda_1, \hdots, \lambda_n$ denote the eigenvalues of $M$. Since $\psi''(s) \leq 0$ for all $s \in \mathbb{R}$, it follows that $\Psi$ is a concave function on the space of symmetric $n \times n$ matrices.

We next rewrite the boundary condition. For each point $x \in \partial \Omega$, we denote by $\nu(x)$ the outward-pointing unit normal vector to $\partial \Omega$ at $x$. Similarly, for each point $y \in \partial \tilde{\Omega}$, we denote by $\tilde{\nu}(y)$ the outward-pointing unit normal vector to $\partial \tilde{\Omega}$ at $y$. By Proposition \ref{uniform.obliqueness}, there exists a positive constant $C_{12}$ such that 
\begin{equation} 
\label{uniform.obliqueness.2}
\langle \nu(x),\tilde{\nu}(\nabla u(x)) \rangle \geq \frac{1}{C_{12}} 
\end{equation}
for all $x \in \partial \Omega$.

We define a subset $\Gamma \subset \partial \Omega \times \mathbb{R}^n$ by 
\[\Gamma = \{(x,y) \in \partial \Omega \times \mathbb{R}^n: \text{$y + t \, \nu(x) \in 
\tilde{\Omega}$ for some $t \in \mathbb{R}$}\}.\] 
For each point $(x,y) \in \Gamma$, we define 
\[\tau(x,y) = \sup \{t \in \mathbb{R}: y + t \, \nu(x) \in \tilde{\Omega}\}\] 
and 
\[\Phi(x,y) = y + \tau(x,y) \, \nu(x) \in \partial \tilde{\Omega}.\] 
If $(x,y)$ lies on the boundary of the set $\Gamma$, then 
\[\langle \nu(x),\tilde{\nu}(\Phi(x,y)) \rangle = 0.\] 
We now define a function $G: \partial \Omega \times \mathbb{R}^n \to \mathbb{R}$ by 
\[G(x,y) = \langle \nu(x),y \rangle - \chi \big ( \langle \nu(x),\tilde{\nu}(\Phi(x,y)) \rangle \big ) 
\, \big [ \langle \nu(x),y \rangle + \tau(x,y) \big ]\] 
for $(x,y) \in \Gamma$ and 
\[G(x,y) = \langle \nu(x),y \rangle\] 
for $(x,y) \notin \Gamma$. Here, $\chi: \mathbb{R} \to [0,1]$ is a smooth cutoff function satisfying $\chi(s) = 1$ for $s \geq \frac{1}{C_{12}}$ and $\chi(s) = 0$ for $s \leq \frac{1}{2 \, C_{12}}$. It is easy to see that $G$ is smooth. Moreover, we have 
\[G(x,y + t \, \nu(x)) = G(x,y) + t\] 
for all $(x,y) \in \partial \Omega \times \mathbb{R}^n$ and all $t \in \mathbb{R}$. Therefore, $G$ is oblique. 

\begin{proposition}
\label{new.problem}
Suppose that $u: \Omega \to \mathbb{R}$ is a convex function such that $\nabla u$ is a diffeomorphism from $\Omega$ to $\tilde{\Omega}$ and $F(D^2 u(x)) = c$ for all $x \in \Omega$. Then $\Psi(D^2 u(x)) = c$ for all $x \in \Omega$. Moreover, we have $G(x,\nabla u(x)) = 0$ for all $x \in \partial \Omega$.
\end{proposition}

\textbf{Proof.} 
It follows from Proposition \ref{interior.c.2.estimate} and Corollary \ref{uniform.convexity} that the eigenvalues of $D^2 u(x)$ lie in the interval $[\frac{1}{C_{11}},C_{10}]$. This implies $\Psi(D^2 u(x)) = F(D^2 u(x)) = c$ for all $x \in \Omega$. 

It remains to show that $G(x,\nabla u(x)) = 0$ for all $x \in \partial \Omega$. In order to verify this, we fix a point $x \in \partial \Omega$, and let $y = \nabla u(x) \in \partial \tilde{\Omega}$. By Proposition \ref{obliqueness}, we have $\langle \nu(x),\tilde{\nu}(y) \rangle > 0$. From this, we deduce that $(x,y) \in \Gamma$ and $\tau(x,y) = 0$. This implies $\Phi(x,y) = y$. Therefore, we have 
\[G(x,y) = \langle \nu(x),y \rangle - \chi \big ( \langle \nu(x),\tilde{\nu}(y) \rangle \big ) \, \langle \nu(x),y \rangle.\] 
On the other hand, it follows from (\ref{uniform.obliqueness.2}) that $\chi(\langle \nu(x),\tilde{\nu}(y) \rangle) = 1$. Thus, we conclude that $G(x,y) = 0$. \\

In view of Proposition \ref{new.problem} we may invoke general regularity results of Lieberman and Trudinger. By Theorem 1.1 in \cite{Lieberman-Trudinger}, the second derivatives of $u$ are uniformly bounded in $C^\gamma(\overline{\Omega})$ for some $\gamma \in (0,1)$. Higher regularity follows from Schauder estimates. \\

\section{The linearized operator}

In this section, we show that all solutions of $(\star)$ are non-degenerate. To prove this, we fix a real number $\gamma \in (0,1)$. Consider the Banach spaces 
\[\mathcal{X} = \bigg \{ u \in C^{2,\gamma}(\overline{\Omega}): \int_\Omega u = 0 \bigg \}\] and 
\[\mathcal{Y} = C^\gamma(\overline{\Omega}) \times C^{1,\gamma}(\partial \Omega).\] 
We define a map $\mathcal{G}: \mathcal{X} \times \mathbb{R} \to \mathcal{Y}$ by 
\[\mathcal{G}(u,c) = \Big ( F(D^2 u) - c, \, (\tilde{h} \circ \nabla u)|_{\partial \Omega} \Big ).\] 
Hence, if $(u,c) \in \mathcal{X} \times \mathbb{R}$ is a solution of $(\star)$, then $\mathcal{G}(u,c) = (0,0)$. \\

\begin{proposition}
\label{linearized.operator.is.invertible}
Suppose that $(u,c) \in \mathcal{X} \times \mathbb{R}$ is a solution to $(\star)$. Then the linearized operator $D\mathcal{G}_{(u,c)}: \mathcal{X} \times \mathbb{R} \to \mathcal{Y}$ is invertible.
\end{proposition}

\textbf{Proof.} 
The linearized operator $\mathcal{B} = D\mathcal{G}_{(u,c)}$ is given by 
\[\mathcal{B}: \mathcal{X} \times \mathbb{R} \to \mathcal{Y}, \quad (w,a) \mapsto (Lw - a,Nw).\] 
Here, the operator $L: C^{2,\gamma}(\overline{\Omega}) \to C^\gamma(\overline{\Omega})$ is defined by 
\[Lw(x) = \text{\rm tr} \Big [ \big ( I + (D^2 u(x))^2 \big )^{-1} \, D^2 w(x) \Big ]\] 
for $x \in \Omega$. Moreover, the operator $N: C^{2,\gamma}(\overline{\Omega}) \to C^{1,\gamma}(\partial \Omega)$ is defined by 
\[Nw(x) = \langle \nabla w(x),\nabla \tilde{h}(\nabla u(x)) \rangle\] 
for $x \in \partial \Omega$. Clearly, $L$ is an elliptic operator. 
Since $u$ is a solution of $(\star)$, Proposition \ref{obliqueness} implies that $\langle \nabla h(x),\nabla \tilde{h}(\nabla u(x)) > 0$ for all $x \in \partial \Omega$. Hence, the boundary condition is oblique. 

We claim that $\mathcal{B}$ is one-to-one. To see this, we consider a pair $(w,a) \in \mathcal{X} \times \mathbb{R}$ such that $\mathcal{B}(w,a) = (0,0)$. This implies $Lw(x) = a$ for all $x \in \Omega$ and $Nw(x) = 0$ for all $x \in \partial \Omega$. Hence, the Hopf boundary point lemma (cf. \cite{Gilbarg-Trudinger}, Lemma 3.4) implies that $w = 0$ and $a = 0$. 

It remains to show that $\mathcal{B}$ is onto. To that end, we consider the operator 
\[\tilde{\mathcal{B}}: \mathcal{X} \times \mathbb{R} \to \mathcal{Y}, \quad 
(w,a) \mapsto (Lw,Nw + w + a).\] It follows from Theorem 6.31 in \cite{Gilbarg-Trudinger} that $\tilde{\mathcal{B}}$ is invertible. Moreover, the operator 
\[\tilde{\mathcal{B}} - \mathcal{B}: \mathcal{X} \times \mathbb{R} \to \mathcal{Y}, \quad (w,a) \mapsto (a,w + a)\] 
is compact. Since $\mathcal{B}$ is one-to-one, it follows from the Fredholm alternative (cf. \cite{Gilbarg-Trudinger}, Theorem 5.3) that $\mathcal{B}$ is onto. This completes the proof. \\

\section{Existence of a solution to $(\star)$} 

In this section, we prove the existence of a solution to $(\star)$. To that end, we employ the continuity method. Let $\Omega$ and $\tilde{\Omega}$ be uniformly convex domains in $\mathbb{R}^n$ with smooth boundary. By Proposition \ref{boundary.defining.function}, we can find a smooth function $h: \Omega \to \mathbb{R}$ with the following properties: 
\begin{itemize}
\item $h$ is uniformly convex
\item $h(x) = 0$ for all $x \in \partial \Omega$
\item If $s$ is sufficiently close to $\inf_\Omega h$, then the sub-level set $\{x \in \Omega: h(x) \leq s\}$ is a ball.
\end{itemize}
Similarly, there exists a smooth function $\tilde{h}: \tilde{\Omega} \to \mathbb{R}$ such that: 
\begin{itemize}
\item $\tilde{h}$ is uniformly convex
\item $\tilde{h}(y) = 0$ for all $y \in \partial \tilde{\Omega}$
\item If $s$ is sufficiently close to $\inf_{\tilde{\Omega}} \tilde{h}$, then the sub-level set $\{y \in \tilde{\Omega}: \tilde{h}(y) \leq s\}$ is a ball.
\end{itemize} 
Without loss of generality, we may assume that $\inf_\Omega h = \inf_{\tilde{\Omega}} \tilde{h} = -1$. For each $t \in (0,1]$, we define 
\[\Omega_t = \{x \in \Omega: h(x) \leq t-1\}, \quad \tilde{\Omega}_t = \{y \in \tilde{\Omega}: \tilde{h}(y) \leq t-1\}.\] 
Note that $\Omega_t$ and $\tilde{\Omega}_t$ are uniformly convex domains in $\mathbb{R}^n$ with smooth boundary. We then consider the following problem (cf. \cite{Brendle}): \\

$(\star_t)$ \textit{Find a convex function $u: \Omega \to \mathbb{R}$ and a constant $c \in (0,\frac{n\pi}{2})$ such that $\nabla u$ is a diffeomorphism from $\Omega_t$ to $\tilde{\Omega}_t$ and $F(D^2 u(x)) = c$ for all $x \in \Omega_t$.} \\

If $t \in [0,1)$ is sufficiently small, then $\Omega_t$ and $\tilde{\Omega}_t$ are balls in $\mathbb{R}^n$. Consequently, $(\star_t)$ is solvable if $t \in (0,1]$ is sufficiently small. In particular, the set 
\[I = \{t \in (0,1]: \text{\rm $(\star_t)$ has at least one solution}\}\] 
is non-empty. It follows from the a-priori estimates in Section 2 that $I$ is a closed subset of $(0,1]$. Moreover, Proposition \ref{linearized.operator.is.invertible} implies that $I$ is an open subset of $(0,1]$. Consequently, $I = (0,1]$. This completes the proof of Theorem \ref{existence}.

\section{Uniqueness}

In this final section, we show that the solution to $(\star)$ is unique up to addition of constants. To that end, we use a trick that we learned from J.~Urbas. 

As above, let $\Omega$ and $\tilde{\Omega}$ be uniformly convex domains in $\mathbb{R}^n$ with smooth boundary. Moreover, suppose that $(u,c)$ and $(\hat{u},\hat{c})$ are solutions to $(\star)$. We claim that the function $\hat{u} - u$ is constant. 

Suppose this is false. Without loss of generality, we may assume that $\hat{c} \geq c$. (Otherwise, we interchange the roles of $u$ and $\hat{u}$.) For each point $x \in \Omega$, we define a symmetric $n \times n$-matrix $B(x) = \{b_{ij}(x): 1 \leq i,j \leq n\}$ by 
\[B(x) = \int_0^1 \Big [ I + \big ( s \, D^2 \hat{u}(x) + (1-s) \, D^2 u(x) \big )^2 \Big ]^{-1} \, ds.\] 
Clearly, $B(x)$ is positive definite for all $x \in \Omega$. Moreover, we have 
\begin{align*} 
&\sum_{i,j=1}^n b_{ij}(x) \, (\partial_i \partial_j \hat{u}(x) - \partial_i \partial_j u(x)) \\ 
&= F(D^2 \hat{u}(x)) - F(D^2 u(x)) = \hat{c} - c \geq 0 
\end{align*}
for all $x \in \Omega$. By the maximum principle, the function $\hat{u} - u$ attains its maximum at a point $x_0 \in \partial \Omega$. By the Hopf boundary point lemma (see \cite{Gilbarg-Trudinger}, Lemma 3.4), there exists a real number $\mu > 0$ such that $\nabla \hat{u}(x_0) - \nabla u(x_0) = \mu \, \nabla h(x_0)$. Using Proposition \ref{obliqueness}, we obtain 
\[\langle \nabla \hat{u}(x_0) - \nabla u(x_0),\nabla \tilde{h}(\nabla u(x_0)) \rangle = \mu \, \langle \nabla h(x_0),\nabla \tilde{h}(\nabla u(x_0)) \rangle > 0.\] 
On the other hand, we have 
\[\langle \nabla \hat{u}(x_0) - \nabla u(x_0),\nabla \tilde{h}(\nabla u(x_0)) \rangle \leq \tilde{h}(\nabla \hat{u}(x_0)) - \tilde{h}(\nabla u(x_0)) = 0\] 
since $\tilde{h}$ is convex. This is a contradiction. Therefore, the function $\hat{u} - u$ is constant. \\

\appendix

\section{The construction of the boundary defining function}

The following result is standard. We include a proof for the convenience of the reader.

\begin{proposition}
\label{boundary.defining.function}
Let $\Omega$ be a uniformly convex domain in $\mathbb{R}^n$ with smooth boundary. Then there 
exists a smooth function $h: \Omega \to \mathbb{R}$ with the following properties:
\begin{itemize}
\item $h$ is uniformly convex
\item $h(x) = 0$ for all $x \in \partial \Omega$
\item If $s$ is sufficiently close to $\inf_\Omega h$, then the sub-level set $\{x \in \Omega: h(x) \leq s\}$ is a ball.
\end{itemize}
\end{proposition} 

\textbf{Proof.} 
Let $x_0$ be an arbitrary point in the interior of $\Omega$. We define a function $h_1: \Omega \to \mathbb{R}$ by 
\[h_1(x) = \frac{d(x,\partial \Omega)^2}{4 \, \text{\rm diam}(\Omega)} - d(x,\partial \Omega).\] 
Since $\Omega$ is uniformly convex, there exists a positive real number $\varepsilon$ such that $h_1$ is smooth and uniformly convex for $d(x,\partial \Omega) < \varepsilon$. We assume that $\varepsilon$ is chosen so that $d(x_0,\partial \Omega) > \varepsilon$. We next define a function $h_2: \Omega \to \mathbb{R}$ by 
\[h_2(x) = \frac{\varepsilon \, d(x_0,x)^2}{4 \, \text{\rm diam}(\Omega)^2} - \frac{\varepsilon}{2}.\] 
For each point $x \in \partial \Omega$, we have $h_1(x) = 0$ and $h_2(x) \leq -\frac{\varepsilon}{4}$. Moreover, for $d(x,\partial \Omega) \geq \varepsilon$, we have $h_1(x) \leq -\frac{3\varepsilon}{4}$ and $h_2(x) \geq -\frac{\varepsilon}{2}$. 

Let $\Phi: \mathbb{R} \to \mathbb{R}$ be a smooth function satisfying $\Phi''(s) \geq 0$ for all $s \in \mathbb{R}$ and $\Phi(s) = |s|$ for $|s| \geq \frac{\varepsilon}{16}$. We define a function $h: \Omega \to \mathbb{R}$ by 
\[h(x) = \frac{h_1(x) + h_2(x)}{2} + \Phi \Big ( \frac{h_1(x) - h_2(x)}{2} \Big ).\] 
If $x$ is sufficiently close to $\partial \Omega$, then we have $h(x) = h_1(x)$. In particular, we have $h(x) = 0$ for all $x \in \partial \Omega$. Moreover, we have $h(x) = h_2(x)$ for $d(x,\partial \Omega) \geq \varepsilon$. Hence, the function $h$ is smooth and uniformly convex for $d(x,\partial \Omega) \geq \varepsilon$.

We claim that the function $h$ is smooth and uniformly convex on all of $\Omega$. To see this, we consider a point $x$ with $d(x,\partial \Omega) < \varepsilon$. The Hessian of $h$ at the point $x$ is given by 
\begin{align*} 
\partial_i \partial_j h(x) 
&= \frac{1}{2} \, \Big [ 1 + \Phi' \Big ( \frac{h_1(x) - h_2(x)}{2} \Big ) \Big ] \, \partial_i \partial_j h_1(x) \\ 
&+ \frac{1}{2} \, \Big [ 1 - \Phi' \Big ( \frac{h_1(x) - h_2(x)}{2} \Big ) \Big ] \, \partial_i \partial_j h_2(x) \\ &+ \frac{1}{4} \, \Phi'' \Big ( \frac{h_1(x) - h_2(x)}{2} \Big ) \, (\partial_i h_1(x) - \partial_i h_2(x)) \, (\partial_j h_1(x) - \partial_j h_2(x)). 
\end{align*}
Note that $|\Phi'(s)| \leq 1$ and $\Phi''(s) \geq 0$ for all $s \in \mathbb{R}$. Since $h_1$ and $h_2$ are uniformly convex, it follows that $h$ is uniformly convex.

It remains to verify the last statement. The function $h$ attains its minimum at the point $x_0$. Therefore, we have $\inf_\Omega h = -\frac{\varepsilon}{2}$. Suppose that $s$ is a real number satisfying 
\[-\frac{\varepsilon}{2} < s < \frac{\varepsilon \, (d(x_0,\partial \Omega) - \varepsilon)^2}{4 \, \text{\rm diam}(\Omega)^2} - \frac{\varepsilon}{2}.\] 
Then we have $\{x \in \Omega: h(x) \leq s\} = \{x \in \Omega: h_2(x) \leq s\}$. Consequently, the set $\{x \in \Omega: h(x) \leq s\}$ is a ball. This completes the proof of Proposition \ref{boundary.defining.function}.

\end{document}